\documentclass[11pt,a4paper]{amsart}

\usepackage{color}

\begin{document}
\newtheorem{cor}{Corollary}[section]
\newtheorem{theorem}[cor]{Theorem}
\newtheorem{prop}[cor]{Proposition}
\newtheorem{lemma}[cor]{Lemma}
\theoremstyle{definition}
\newtheorem{defi}[cor]{Definition}
\newtheorem{remark}[cor]{Remark}
\newtheorem{example}[cor]{Example}

\newcommand{\cC}{\mathcal{C}}
\newcommand{\cH}{\mathcal{H}}
\newcommand{\cun}{\cC^{\infty}}
\newcommand{\cz}{{\mathbb C}}
\newcommand{\fm}{f^{-1}}
\renewcommand{\index}{\mathrm{index}}
\newcommand{\oX}{\bar{X}}
\newcommand{\oD}{\overline{D}}
\newcommand{\px}{\partial_x}
\newcommand{\rz}{{\mathbb R}}
\newcommand{\Spec}{\operatorname{Spec}}
\newcommand{\supp}{\mathrm{supp}}
\newcommand{\tf}{{\tilde{f}}}
\newcommand{\tg}{{\tilde{g}}}
\newcommand{\tih}{\tilde{h}}
\newcommand{\tphi}{\tilde{\phi}}
\newcommand{\vol}{\mathrm{vol}}
\newcommand{\Xe}{X_\epsilon}

\def\f{\varphi}
\def\e{\varepsilon}
\def\d{\delta}

\title[Dirac spectrum and gradient conformal vector fields]{The Dirac
spectrum on manifolds with gradient conformal
vector fields}
\author{Andrei Moroianu}
\address{Centre de Math\'ematiques\\
\'Ecole Polytechnique\\
91128 Palaiseau Cedex\\
France}
\email{am@math.polytechnique.fr}
\author{Sergiu Moroianu}
\address{Institutul de Matematic\u{a} al Academiei Rom\^{a}ne\\
P.O. Box 1-764\\RO-014700
Bu\-cha\-rest, Romania}
\address{\c Scoala Normal\u a Superioar\u a 
Bucharest, calea Grivi\c tei 21, Bucharest, Romania}
\email{moroianu@alum.mit.edu}
\thanks{The second author was partially supported 
from the contracts 2-CEx06-11-18/2006 and CNCSIS-GR202/19.09.2006.}
\date{\today}
\begin{abstract}
We show that the Dirac operator on a spin manifold 
does not admit $L^2$ eigenspinors
provided the metric has a certain
asymptotic behaviour and is a warped product near infinity.
These conditions on the metric are fulfilled in particular 
if the manifold is complete and carries a non-complete
vector field which outside a compact set is gradient conformal and
non-vanishing.
\end{abstract}
\maketitle

\section{Introduction}

The Dirac operator on a closed spin manifold is 
essentially self-adjoint as an unbounded operator in $L^2$, and has purely 
discrete spectrum. Its eigenvalues grow at a certain speed determined by the 
volume of the manifold and its dimension. Hence, although determining the 
eigenvalues can be a daunting task, the nature of the spectrum is rather 
well understood.
 
On non-compact manifolds, the spectrum of the Dirac operator can behave in a 
variety of ways. For instance, the Dirac operator on $\rz^n$ has purely 
absolutely continuous spectrum, so in particular there are no 
$L^2$ eigenspinors. In contrast, B\"ar \cite{baer} showed that on complete 
spin hyperbolic manifolds of finite volume, the spectrum is purely discrete 
if the induced Dirac operators on the ends are invertible. 
In this situation even the classical Weyl law for the distribution of 
the eigenvalues holds \cite{wlom}.
Otherwise, if the limiting Dirac operator is not invertible, then the 
essential spectrum is the whole real axis. Similar results appear in 
\cite{mso} for the Laplace operator on forms and for magnetic Schr\"odinger operators.

In this paper we show that for a class of -- possibly
incomplete -- spin Riemannian manifolds $(X,g)$ which includes certain
hyperbolic   
manifolds, the Dirac operator $D$ does not carry $L^2$ eigenspinors of
real eigenvalue.
In particular, we deduce that the $L^2$ index of the Dirac operator 
on $(X,g)$ vanishes.

Our main result (Theorem \ref{th1}) makes special assumptions on the
metric of $X$. Geometrically, these assumptions imply the existence of a
non-complete vector field on $X$ which is gradient conformal on an
open subset $U$ of $X$. Conversely, we show in Section \ref{gcvf} that 
the existence of such vector fields implies the hypothesis of Theorem
\ref{th1} provided that $X$ is complete and $X\setminus U$ is compact.
As a corollary, we obtain that on a complete spin manifold $X$ which
carries a non-complete vector field which is gradient conformal
outside a compact subset of $X$, the Dirac operator has purely
continuous spectrum (Theorem \ref{main}).

\section{The main result}
Let $(\oX^n,h)$ be a connected spin manifold with boundary
with interior $X$. Assume that
there exists a boundary component $M$ so that $h$ is a product
in a neighbourhood of $M$. Denote by $x:\oX\to[0,\infty)$ the distance
to $M$, so $h=dx^2+h_M$ near $M$. Note that $M$ inherits a spin structure from
$\oX$.

Let $f:X\to(0,\infty)$ be a smooth conformal factor which 
depends only on $x$ in a neighbourhood of $M$. 

\begin{theorem}\label{th1}
Assume that $M$ is at infinite distance from $X$
with respect to the conformal metric 
\begin{equation}\label{met}
g:={f^2}{h}={f(x)^2}({dx^2+h_M}).
\end{equation}
Moreover assume that the Dirac operator $D_M$ on $(M, h_M)$  is 
essentially self-adjoint. Then the Dirac operator of $(X,g)$
does not have any distributional $L^2$ eigenspinors of real eigenvalue.
\end{theorem}

Lott \cite{lott} proved that there is no $L^2$ harmonic spinor 
under somewhat similar assumptions. Namely,
$h$ could be any metric smooth up to $M$, and $f$ 
could vary in the $M$ directions. 
However, Lott assumes that $f^{-1}$ extends to a locally Lipschitz function 
on $\oX$ which eventually must be locally bounded by a multiple of $x$, 
while our hypothesis only asks that
\begin{equation}\label{intf}
\int_0^\epsilon {f(x)}{dx}=\infty.\end{equation}
In particular, unlike in \cite{lott}, the function $f^{-1}$ may be
unbounded near $M$. 

If we assume that $(X,g)$ is complete as in \cite{lott}, we deduce that
$(M,h_M)$ is complete so $D_M$ is essentially self-adjoint. With this
assumption 
we also know that $D_g$ is essentially self-adjoint, so its spectrum is real.
But we do not need to make this assumption as the statement is 
``local near $M$''. This seems to be also the case in \cite{lott},
although it is not claimed explicitly. However, there are many instances 
of \emph{incomplete} manifolds whose Dirac operator is essentially 
self-adjoint, cf.\ Example \ref{exsa}.

The rest of this section is devoted to the proof of Theorem
\ref{th1}.

It is known since Hitchin \cite{hitchin} that the Dirac operator 
has a certain conformal invariance property. More precisely, if $g=f^2h$, 
the Dirac operator $D_g$ is conjugated to $f^{-1}D_h$ by the Hilbert 
space isometry
\begin{align}\label{cocha}
L^2(X,\Sigma,\vol_g)&\to L^2(X,\Sigma, f\vol_h),& \psi&\mapsto 
f^{\frac{n-1}{2}}\psi,
\end{align}
(see \cite[Proposition 1]{nistor}). Let $f^{-\frac{n-1}{2}}\phi
\in L^2(X,\Sigma,\vol_g)$ be an eigenspinor of $D_g$ (in the sense of
distributions) of eigenvalue $l\in\rz$. Then 
\begin{equation}\label{eqei}
(f^{-1}D_h-l)\phi=0.
\end{equation}
By elliptic regularity, $\phi$ is in fact a smooth spinor on $X$ so
\eqref{eqei} is equivalent to 
\begin{equation}\label{eqd}
D_h\phi=lf\phi.
\end{equation}

Let $c^0$ denote the Clifford multiplication by the unit normal vector
$\px$ with respect to $h$. Then $D_h$ decomposes near $M$ as follows:
\[D_h=c^0(\px+A)\]
where for each $x>0$, $A$ is a differential 
operator on the sections of $\Sigma$ over $M\times \{x\}$; moreover $A$
is independent of $x$. Note that $\Sigma\big|_{M}$ is either the spinor bundle
$\Sigma(M)$ of $M$ with respect to the induced spin structure 
(if $n$ is odd), or two copies of $\Sigma(M)$ if $n$ is even.
We can describe $A$ in terms of the Dirac operator $D_M$ on
$M$ with respect to the metric $h_M$ as follows:
\begin{equation}\label{fA}
A=\begin{cases}
 D_M &\text{ for $n$ odd,}\\
\begin{bmatrix} D_M&0\\0&-D_M\end{bmatrix}&\text{for $n$ even}.
\end{cases}
\end{equation}
In both cases, $A$ is symmetric and elliptic. Since $D_M$ is essentially
self-adjoint, so is $A$ and we use the same symbol for its unique self-adjoint
extension.

\subsection*{The case where $A$ has pure-point spectrum}
For the sake of clarity we make temporarily the assumption that 
$L^2(M,\Sigma,\vol_{h_M})$ admits an orthonormal basis  made of eigenspinors 
of $A$ of real eigenvalue (equivalently, $A$ has pure-point spectrum). 
Since $A$ is essentially $1$ or $2$ copies of $D_M$, this happens 
if $M$ is compact, but also more generally.
\begin{example} \label{exsa}
Suppose that either $(M,h_M)$ is 
conformal to a closed cusp metric such that the induced Dirac operators 
on the ends are invertible (see \cite{wlom}), or that
$(M,h_M)$ is compact with isolated conical singularities and the 
Dirac eigenvalues 
on the cone section do not belong to $(-\frac12,\frac12)$ (see {\em e.g.}
\cite{chou}). Then the operator $D_M$ is essentially self-adjoint 
with purely discrete spectrum.
\end{example}

Note that since $c^0$
and $A$ anti-commute, we have 
\[A\phi_\lambda=\lambda\phi_\lambda \Longrightarrow
Ac^0\phi_\lambda=-\lambda c^0\phi_\lambda\]
and so the spectrum of $A$ is symmetric around $0$. 
For each $x$, decompose $\phi$ onto the positive 
eigenspaces of $A$ as follows:
\[
\phi=\sum_{\substack{\lambda\in\Spec A\\
\lambda>0}} \left(a_\lambda(x)\phi_\lambda +b_\lambda(x)c^0\phi_\lambda\right)
+\phi_0(x)\]
where $\phi_0(x)\in\ker A$ for all $x>0$, $\lambda>0$ is an eigenvalue of $A$,
and
$\phi_\lambda$ is an eigenspinor of eigenvalue $\lambda$ and norm $1$
in $L^2(M,\Sigma,\vol_{h_M})$. We compute
\[\begin{split}
\|\phi\|^2_{L^2(X,\Sigma,{f} \vol_{h})}&=
\int_{\{x>\epsilon\}} |\phi|^2 {f} \vol_h\\
&\quad +\int_0^\epsilon\int_M |\phi_0(x)|^2 \vol_{h_M} {f}(x)dx\\
&\quad+\sum_{\substack{\lambda\in\Spec A\\
\lambda>0}}\int_0^\epsilon \left(|a_\lambda(x)|^2 + |b_\lambda(x)|^2\right)
{f}(x) dx.
\end{split}\]
In particular, $\phi\in L^2(X,\Sigma,{f} \vol_h)$ implies that
\begin{align}
\int_0^\epsilon\int_M |\phi_0(x)|^2 \vol_{h_M} {f}(x)dx<\infty\label{pold}\\
\intertext{and}
a_\lambda, b_\lambda\in L^2((0,\epsilon), {f} dx). \label{alld}
\end{align}
The eigenspinor equation \eqref{eqd} becomes
\[\begin{split}
0&=(D_h-l{f})\phi\\&=c^0(\px+A+c^0l{f})\phi\\
&=c_0\sum_{\substack{\lambda\in\Spec A\\ \lambda>0}}
\left((a_\lambda '+\lambda a_\lambda -l{f} b_\lambda)\phi_\lambda
+(b_\lambda '-\lambda b +l{f} a_\lambda)c^0\phi_\lambda\right)\\
&\quad+c^0\phi_0 '-l{f}\phi_0.
\end{split}\]
So we get
\begin{equation}\label{pk}
\phi_0 '=-l{f} c^0\phi_0
\end{equation}
for the part of $\phi$ which for all fixed $x$ lives in $\ker A$, while for 
$\lambda>0$,
\begin{equation}\label{ldz}
\begin{cases}a_\lambda '=-\lambda a_\lambda+l{f} b_\lambda\\
b_\lambda '=-l{f} a_\lambda+\lambda b_\lambda.\end{cases}
\end{equation}
First, since $l$ is real and $c^0$ is skew-adjoint, Eq.\ \eqref{pk}
implies that the function $|\phi_0(x)|^2$ is constant in $x$. Together with 
\eqref{intf} and \eqref{pold}, we see that $\int_M |\phi_0(x_0)|^2dh_M=0$
so $\phi_0\equiv 0$.

We show now that for all $\lambda>0$, 
the system \eqref{ldz} does not have nonzero solutions satisfying
\eqref{alld}, {\em i.e.}, in $L^2((0,\epsilon),{f} dx)$. 
\begin{remark}
The Wronskian of \eqref{ldz} is constant in $x$, so by \eqref{intf},
the two fundamental solutions cannot belong  simultaneously to
$L^2((0,\epsilon),{f} dx)$. However, this
fact alone does not stop \emph{one} solution from being in $L^2$!
\end{remark}
Fix $0<\lambda\in\Spec A$ and set 
\begin{align*}
a(x)&:=e^{\lambda x} a_\lambda(x),&
b(x)&:=e^{-\lambda x}b_\lambda(x).
\end{align*}
Then \eqref{ldz} becomes
\begin{equation}\label{ldzp}
\begin{cases}
a'(x)=le^{2\lambda x} {f}(x) b(x)\\
b'(x)=-le^{-2\lambda x}{f}(x) a(x).\end{cases}
\end{equation}
Note that the system \eqref{ldzp} has real coefficients (here we use 
the hypothesis that $l$ is real) so by splitting into real and 
imaginary parts, we can assume that $a,b$ are also real.

Since $e^{\pm 2\lambda x}$ is bounded for $0\leq x\leq \epsilon<\infty$,
condition \eqref{alld} implies that $a,b\in L^2((0,\epsilon),{f} dx)$.
If $l=0$ then $a,b$ are constant functions, which by \eqref{intf}
do not belong to
$L^2((0,\epsilon),{f} dx)$ unless they are $0$. So in that case $a_\lambda$
and $b_\lambda$ vanish identically. 

Let $L^1$ denote the space of integrable functions on $(0,\epsilon)$
with respect to the Lebesgue measure $dx$. Then \eqref{alld}
implies that ${f} ab\in L^1$. So from \eqref{ldzp},
\[(a^2)'=2le^{2\lambda x} {f}(x)a(x) b(x)\in L^1.\]
Hence 
\[\lim_{x\to 0} a^2(x)=a^2(x_0)-\int_0^{x_0} (a^2)'(t) dt\]
exists, in other words $a^2$ extends continuously in $x=0$. The same argument
shows that $b^2$ (and so also $a,b$) are continuous in $0$.

The case $l=0$ was treated above so we can assume that $l\neq 0$.
We claim that $b(0)=0$. Otherwise, by continuity, $b(x)\neq 0$ for 
$0\leq x\leq x_0$
so from \eqref{ldzp}, there exists $C>0$ with
$|a'(x)|>C{f}(x)$ (in particular, by continuity $a'$ has constant
sign) and therefore  
\[|a(0)-a(x_0)|=\lim_{x\to 0}\int_x^{x_0} |a'(x)|dx>C\int_0^{x_0}
{f}(t)dt=\infty\] 
which is a contradiction.  So $b(0)=0$ and similarly $a(0)=0$.

We pull now our final trick. Recall that $\lambda >0$. Consider the function 
\[F(x):=e^{-4\lambda x}a^2(x)+b^2(x)\geq 0.\]
We just showed that $F$ is continuous in $0$
and $F(0)=0$. From \eqref{ldzp} we compute 
\[F'(x)=-4\lambda e^{-4\lambda x} a^2(x)\leq 0.\]
By collecting what we know about $F$, we note:
\begin{enumerate}
\item $F(x)\geq 0$;
\item $F(0)=0$;
\item $F'(x)\leq 0$ for $\epsilon>x>0$
\end{enumerate} 
Together these facts imply that $F\equiv 0$ on $(0,\epsilon)$.
This is equivalent to saying $a(x)=b(x)=0$ for all $0\leq x\leq \epsilon$.

So we showed that $\phi$ vanishes near $M$.
The eigensections of $D$ have the unique continuation property, which 
implies that $\phi$ vanishes on $X$.

\begin{remark}
Recall that the spectrum of $A$ is symmetric around $0$, and 
$a_\lambda,b_\lambda$ are the coefficients in $\phi$ of the eigensections
$\phi_\lambda$, $c^0\phi_\lambda$ of eigenvalue 
$\lambda$, respectively $-\lambda$.
It may seem that starting the decomposition using positive $\lambda$ was a 
fortunate choice, otherwise the last argument would not hold. But in fact, 
the argument works for $-\lambda$ by choosing a 
different function $\tilde{F}(x):=a^2(x)+e^{4\lambda x}b^2(x)$.
\end{remark}

\subsection*{The general case}
Let us remove the  assumption that the spectrum of $A$ is 
purely discrete.
We will model the proof on the argument given above, which is now loaded with
technical subtleties.

By assumption, a neighbourhood of the infinity in $X$ 
is isometric to $(0,\epsilon)\times M$ with the metric \eqref{met}. 

After the unitary transformation \eqref{cocha},
the eigenspinor equation reads as before
\begin{equation}\label{egsp}
(\px+A)\phi=-c^0f^{-1}l\phi
\end{equation}
where $l$ is real, and $\phi$ is smooth (by elliptic regularity) 
and square-integrable. To make this last condition precise, let
$I_\epsilon$ denote the interval $(0,\epsilon)$ with the measure
$f(x)dx$. 
Denote by $\cH$ the Hilbert space $L^2(M,\Sigma,\vol_{h_M})$, then
\begin{equation}\label{ld}
\phi\in L^2(\Xe,\Sigma, {f}\vol_h)=L^2(I_\epsilon,\cH).
\end{equation}
In particular for almost all $x\in I_\epsilon$, we have $\phi_x\in\cH$.

Let $\chi$ be the characteristic function of the interval $[-N,N]$ for some
$N\in\rz$. Let $\chi(A)$ be the corresponding spectral projection. Since 
$A$ anti-commutes with $c^0$ and $\chi$
is even, it follows that $\chi(A)$ commutes with $c^0$.

Let $\cH^1\subset \cH$  be the domain of $A$ and $\cH^{-1}\supset\cH$ 
its dual inside distributions, {\em i.e.}, the space of those distributions 
which extend continuously to $\cH^1$. Since $\chi$ has compact support, 
we deduce that $\chi(A)$ acts continuously from $\cH$ to $\cH^1$, and also
from $\cH^{-1}$ to $\cH$. 

From \eqref{ld} we deduce $A\phi\in L^2(I_\epsilon,\cH^{-1})$,
$\chi(A)\phi\in L^2(I_\epsilon,\cH^{1})$ and
\[\chi(A)(A\phi)=A(\chi(A)\phi)\in L^2(I_\epsilon,\cH).\]
Similarly, 
$\px\phi\in H^{-1}_{loc}(I_\epsilon,\cH)$ and 
\[\chi(A)(\px\phi)=\px(\chi(A)\phi)\in H^{-1}_{loc}(I_\epsilon,\cH^1).\]
It follows that $\tphi:=\chi(A)\phi$
satisfies (in distributions) the eigenspinor equation \eqref{egsp}. 
Denote by $\cH_N$ the range of the projection $\chi(A)$, then
$\tphi\in L^2(I_\epsilon,\cH_N)$. Most importantly for us, $A$
acts as a self-adjoint \emph{bounded} operator on $\cH_N$.

\begin{lemma}\label{lc}
Let $H$ be a Hilbert space, $A:H\to H$ a bounded self-adjoint
operator, and $c^0$ a skew-adjoint 
involution of $H$ which anti-commutes with $A$. Then for every $l\in\rz$, 
the equation $(\ref{egsp})$ does not have $($distributional$)$
solutions in $L^2(I_\epsilon,H)$ other than 0.
\end{lemma}
\begin{proof} Let $\phi$ be a solution of (\ref{egsp}), square-integrable
with respect to the measure ${f(x)} dx$ on $I_\epsilon$. By elliptic
regularity,  
$\phi$ is smooth in $x$.
Since $\exp(xA)c^0=c^0\exp(-xA)$, we get $\px(\exp(xA)\phi)=
-l{f} c^0\exp(-xA)\phi$, hence the family of $H$-norms
$x\mapsto \|\px(\exp(xA)\phi_x)\|$ is square-integrable with respect to
the measure $f^{-1}dx$. Since $x\mapsto\|\phi_x\|$ is 
$L^2$ with respect to ${f} dx$, and $\exp(xA)$ is uniformly bounded, 
by the Cauchy-Schwartz inequality we see
that the function 
\[x\mapsto \frac{d}{dx}\|\exp(xA)\phi_x\|^2\] is integrable with respect to 
the Lebesgue measure $dx$. 
Thus 
\[\|\exp(xA)\phi_x\|^2=\|\exp(x_0A)\phi_{x_0}\|^2+
\int_{x_0}^{x}\frac{d}{dx}\|\exp(xA)\phi_x\|^2 dx\]
has a finite limit as $x\searrow 0$. We claim that this limit is
$0$. Otherwise, since $\lim_{x\searrow 0} \exp(xA)=1$, we would have
$\lim_{x\searrow 0} \|\phi_x\|^2>0$ which, together with \eqref{intf},
contradicts the fact that $\phi$ is square-integrable with respect to ${f} dx$.

Thus $\phi_x$ tends in norm to $0$ in $H$ as $x\searrow 0$. Let now
$|A|$ be the absolute value of $A$, and define
\[F(x):=\| \exp(-x|A|) \phi_x)\|^2.\]
We notice that $c^0$ commutes with $|A|$ since it commutes with
$A^2$. A direct computation shows, using that $c^0$ is skew-adjoint,
\[\frac{dF}{dx}=-2\langle(A+|A|)\exp(-x|A|) \phi_x,\exp(-x|A|) 
\phi_x\rangle\leq 0.\]
Hence $F$ is decreasing, on the other hand it vanishes at $x=0$ and it is 
non-negative, so in conclusion it vanishes identically. Since 
$\exp(-x|A|)$ is invertible, we conclude that $\phi\equiv 0$.
\end{proof}
We apply this lemma to the eigenspinor $\tphi$ constructed 
above with $H=\cH_N$. Therefore $\chi(A)\phi=0$ for all $N\in\rz$. But
as $N\to\infty$ we have $\chi(A)\phi\to\phi$. By the uniqueness of the limit,
$\phi$ must be identically zero on $(0,\epsilon)\times M$ which is an 
open subset of $X$. By the unique continuation property, it follows that 
$\phi$ vanishes on $X$ as claimed.
This ends the proof of Theorem \ref{th1}.
\begin{remark}
Theorem \ref{th1} holds
for Dirac operators twisted with a Hermitian bundle $E$ whose connection 
$\nabla^E$ is flat in the direction of the conformal gradient vector 
field $\xi=\px$, {\em i.e.}, such that the contraction of 
the curvature of $\nabla^E$ with the field $\xi$ vanishes. 
We only need to replace in equation \eqref{fA} the operator $D_M$ 
by the twisted operator $D_M^E$. The flatness condition ensures that 
this operator is independent of $x$. The rest of the proof remains unchanged.
\end{remark}

\section{A formal extension of Theorem \ref{th1}}

For applications, it might be useful to view the metric $g$ given
by (\ref{g}) in different coordinates. We state below the most
general reformulation of Theorem \ref{th1}.

\begin{theorem}\label{cor2}
Let $\oX$ be a smooth manifold with boundary, let $M$
be a boundary component,  
and let $[0,\epsilon)\times M\hookrightarrow \oX$ be a collar
neighbourhood of $M$.  
Take a smooth $L^1$ function $\rho:(0,\epsilon)\to(0,\infty)$
and let $\tih$ be a possibly incomplete Riemannian metric on $X$ 
which is a warped product near $M$: 
\[\tih=dx^2+\rho^{-2}(x)h_M.\]
Let $f:X\to(0,\infty)$ be a smooth conformal factor 
depending only on $x$ near $M$ and satisfying $\int_0^\epsilon
f(x)dx=\infty$. Assume that the Dirac operator $D_M$ on $(M,h_M)$ 
with the induced spin structure is essentially self-adjoint. 
Then the Dirac operator $D_\tg$ of the metric 
$\tg:={f^2}{\tih}$
does not carry square-integrable eigenspinors of real eigenvalue.
\end{theorem}
\begin{proof}
The metric $\tih$ is conformal to $(\rho(x)dx)^2+h_M$.
Set 
\[t(x):=\int_0^x \rho(s)ds.\]
Since $\rho$ is in $L^1$, it follows that $t$ is well-defined and $t(0)=0$.
Since $\rho$ is positive, $x\mapsto t(x)$ is an increasing diffeomorphism from
$(0,\epsilon)$ to $(0,t(\epsilon))$ which extends to a homeomorphism
between $[0,\epsilon)$ and $[0,t(\epsilon))$. We write $x=x(t)$ for
its inverse. 
Define $\tf(t):=\frac{f(x(t))}{\rho(x(t))}$.
Clearly, $dt=\rho dx$ so
\[\tg={\tf^2(t)}\left(dt^2 +h_M\right).\]
Note that 
\[\int_0^{t(\epsilon)}{\tf(t)}{dt}=\int_0^\epsilon{f(x)}{dx}=\infty\] 
so we can apply Theorem \ref{th1}.
\end{proof}

\section{Gradient conformal vector fields}
\label{gcvf}

Let $(X^n,g)$ be a Riemannian manifold. A {\em gradient conformal 
vector field} (or GCVF) on $X$ is a conformal vector field $\xi$ which is 
at the same time the gradient of a function on $X$:
\begin{equation}\label{sys}\begin{cases} \mathcal L_\xi g=\alpha g,
\qquad\hbox{for some }\alpha\in{\mathcal C}^\infty (X), \\ \xi=\nabla^g F,
\qquad\hbox{for some }F\in{\mathcal C}^\infty (X).
\end{cases}\end{equation}
Gradient conformal 
vector fields were studied intensively in the 70s
(see \cite{bo} and references therein). More recently, they turned out
to be a very useful tool in understanding other geometric objects,
like closed twistor 2-forms on compact Riemannian manifolds \cite{an}.
The aim of this section is to prove the following:
\begin{theorem}\label{main}
Let $(X,g)$ be a complete spin manifold. Assume that $X$ carries a 
non-complete vector field which outside some compact subset 
is nowhere-vanishing and GCVF. Then the Dirac operator of $(X,g)$
does not carry square-integrable eigenspinors.
\end{theorem}
This theorem is a direct consequence of Theorem \ref{th1} and
Proposition \ref{prop} below.

We first recall some basic properties of GCVFs.

\begin{lemma}\label{lem} Let $\xi$ be a GCVF satisfying the system
  $(\ref{sys})$. Then the following assertions hold:

$(i)$ The covariant derivative of $\xi$ depends on only the co-differential
of $\xi$:
\begin{equation}\label{g}
\nabla_Y\xi=\phi Y,\qquad\forall\ Y\in TX,
\end{equation}
where the function $\phi$ equals $-\frac1n\delta\xi$.

$(ii)$ Let $X_0$ be the set of points where $\xi$ does not vanish. The
distribution $\xi ^\perp$ defined on $X_0$ is involutive and its
maximal integral leaves are exactly the connected components of the
level sets of $F$ on $X_0$. 

$(iii)$ The length of $\xi$ is constant along the integral leaves of
$\xi ^\perp$. 

$(iv)$ The integral curves of $\xi$ are geodesics and $F$ is strictly
increasing along them.

$(v)$ Each point $p$ of $X_0$ has a neighbourhood isometric to
$$((-\e,\e)\times V,f^2(x)(dx^2+h)),$$ 
where $(V,h)$ is a local integral
leaf of $\xi ^\perp$ through $p$ and $f:(-\e,\e)\to\rz^+$ is some
positive function. In these coordinates $\xi$ corresponds to
$\partial/\partial x$ and $f(x)$ is the norm of $\xi$ on the leaf 
$\{x\}\times V$. 
\end{lemma}

\begin{proof} $(i)$ The first equation in (\ref{sys}) is
equivalent to the vanishing of 
the trace-free symmetric part of $\nabla\xi$. The second
equation of (\ref{sys}) implies that the skew-symmetric part 
of $\nabla\xi$ vanishes
too. We are left with $\nabla_X\xi=\phi X$ for some function $\phi$. 
Taking the scalar product with $X$ and the sum over an orthonormal
basis $X=e_i$ yields $\phi=-\frac1n\delta\xi$.

$(ii)$ By definition, the distribution $\xi ^\perp$ on $X_0$ is 
exactly the kernel of the 1-form
$dF$, so it is involutive. Let $M$ be a
maximal integral leaf of $\xi ^\perp$. Clearly $F$ is constant on $M$, 
which is
connected, so $M$ is a subset of some level
set $F^{-1}(y)$. Moreover $M$ is open in $F^{-1}(y)$, as can be seen in
local charts. Thus $M$ is a connected component of $F^{-1}(y)$.
$(iii)$ For every $Y\in TM$ one can write
$$Y(|\xi|^2)=2g(\nabla_Y\xi,\xi)\stackrel{(\ref{g})}{=}\phi
g(Y,\xi)=0,$$
so $|\xi|^2$ is constant on $M$. 

$(iv)$ Let $\f_t$ denote the local flow of $\xi$ and let
$\gamma_t:=\f_t(p)$ for some $p\in M$. Taking $X=\xi$ in (\ref{g})
yields 
$$\nabla_{\dot\gamma}\dot\gamma=\nabla_\xi\xi=\phi\xi=\phi\dot\gamma,$$ 
which 
shows that $\gamma_t$ is a (non-parametrized)
geodesic. Furthermore,
$$\frac{d}{dt}F(\gamma_t)=\dot\gamma_t(F)=\xi(F)=|\xi ^2|>0,$$
so $F(\gamma_t)$ is increasing.

$(v)$ The tangent bundle of $X_0$ has two involutive 
orthogonal distributions
$\rz\xi$ and $\xi^\perp$. The Frobenius integrability theorem 
shows that there exists a local coordinate system $(x,y_1,\ldots,y_{n-1})$
around every $p\in X_0$ such that $\xi=\px$ and $\partial_{y_i}$ span $\xi^\perp$. 
Let $V$ denote the set $\{x=0\}$ in these coordinates. 
The metric tensor can be written
$$g=f^2dx^2+\sum_{i,j=1}^{n-1}g_{ij}dy_i\otimes dy_j.$$
From $(iii)$ we see that $f$ only depends on $x$. Using the first
equation in the system (\ref{sys}) we get 
$$2g_{ij}(\log f)'(x)=\frac{\partial g_{ij}}{\partial
  x},\qquad\forall\ i,j\le n-1,$$ 
which shows that $g_{ij}(x,y)=f^2(x)h_{ij}(y)$ for some metric tensor 
$h$ on $V$.
\end{proof}

It turns out that under some completeness assumptions,
the last statement of the lemma also holds globally:

\begin{prop}\label{prop} Let $(X^n,g)$ be a complete Riemannian
  manifold. If $\xi$ is 
  a non-complete vector field on $X$ which, outside a compact subset of $X$, 
  is gradient conformal and non-vanishing,
  then there exists an open
  subset of $X$ which is isometric to $((0,c)\times M,f^2(x)(dx^2+h))$
  for some complete Riemannian manifold $(M^{n-1},h)$ and smooth
  positive function $f:(0,c)\to\rz^+$ with $\int_0^c f(x)dx=\infty.$ 
\end{prop}
\begin{proof}
Let $\f_t$ denote the local flow of $\xi$ and let $K$ be a compact
subset of $X$ such that $\xi$ is gradient 
conformal and nowhere-vanishing on $X\setminus K$. By definition,
$\xi=\nabla F$ for some function $F$ defined on $X\setminus K$.
Consider the open set
$$K_\e:=\{p\in X\ |\ d(p,K)<\e\}.$$ 
Since $\overline{K_\e}$ is compact, there exists some $\d>0$
such that $\f_t$ is defined on $K_\e$ for every $|t|<\d$. Since $\xi$
is non-complete, there exists some $p\in X$ and $a\in\rz$ such that 
$\f_t(p)$ tends to infinity as $t$ tends to $a$. By changing $\xi$ to
$-\xi$ if necessary, we can assume that $a>0$. From
the definition of $\d$ we see that $\f_t(p)\in
X\setminus K_\e$ for all $t\in [a-\d,a)$. Since 
$$\lim_{t\to a}\int_{0}^t |\xi_{\f_s(p)}|ds\ge \lim_{t\to
  a}d(p,\f_t(p))=\infty,$$
the norm of $\xi$ has to be unbounded along its integral curve
through $p$ in the positive direction. Therefore one can find a point
$q:=\f_{t_0}(p)$ ($t_0\in [0,a)$) on this integral curve such that
$|\xi_q|$ is larger than the supremum of the norm of $\xi$ over
$\overline{K_\e}$. Let $M$ be the maximal leaf through $q$ of the involutive
distribution $\xi ^\perp$ (defined on $X\setminus K$). Since the
norm of $\xi$ is constant on $M$, it is clear 
that $M$ does not intersect $K_\e$. We notice that $M$ is
complete with respect to the induced Riemannian metric $h$. This does not
follow directly from the completeness of $(X,g)$ 
since the distribution $\xi ^\perp$ is only defined and involutive on
$X\setminus K$. Nevertheless, since $M$ is a connected component of
some level set of $F$, it is closed in $X$, and every closed
submanifold of a complete Riemannian manifold is also complete with
respect to the induced Riemannian metric.

From the definition of $q$, it is clear that the integral curve
$\f_t(q)$ is defined for $t<a-t_0$. From Lemma \ref{lem}, 
two integral curves of $\xi$ which do not meet $K$,
which are issued from points of the same maximal leaf, 
are geodesics and have the
same length. Consequently, for every other
point $q'\in M$, the integral curve of $\xi$ in the positive direction
is defined at least for all $t<a-t_0$. Moreover, the map 
$$\psi:M\times (0,a-t_0)\to X ,\qquad \psi(r,t):=\f_t(r)$$
is one-to-one since the vector field $\xi$ 
does not have zeros on $M$.

Finally, Lemma \ref{lem} $(v)$ shows that $\psi$ is
an isometric embedding of $(M\times (0,a-t_0),f^2(x)(dx^2+h))$ into
$(M,g)$, where $f(x)$ denotes the length of $\xi$ on the maximal leaf 
$\f_x(M)$ of $\xi ^\perp$. \end{proof}

\begin{remark}
The incompleteness condition on $\xi$ in Theorem
\ref{main} is necessary. Indeed, complete
hyperbolic manifolds of finite volume are isometric
outside a compact set to a disjoint union of
\emph{cusps}, {\em i.e.}\ cylinders $(0,\infty)\times T$ over
some flat connected Riemannian manifold $(T,h)$, with
metric
\[dt^2+e^{-2t}h=e^{-2t}((d e^{t})^2+h).\]
The vector field $e^{-t}\partial/\partial t$ is GCVF
and complete. These manifolds are known to have
\emph{purely discrete spectrum} if the spin structure
on each cusp is non-trivial \cite{baer}, which is the
case for instance in dimension $2$ or $3$ when there
is only one cusp. The eigenvalues then obey the Weyl
asymptotic law \cite{wlom}. On the contrary, when some
cusps have non-trivial spin structures, the spectrum of
$D$ is the real line. In this case (like for the
scalar Laplace operator) the existence of $L^2$ eigenspinors
is generally unknown.
\end{remark}

\section{Applications}

\subsection*{Real hyperbolic space}
The Poincar\'e disk model of the hyperbolic space is conformally equivalent to
the standard flat metric. In polar coordinates, this metric is a warped product
so Theorem \ref{cor2} shows that the Dirac operator on the hyperbolic space 
does not have point spectrum (the spectrum is real since ${\mathbb H}$
is complete). This was first studied with different methods
by Bunke \cite{bunk}.

\subsection*{Hyperbolic manifolds}
More generally, let $(M^n,h_M)$ be a spin hyperbolic manifold whose 
Dirac operator is essentially self-adjoint. Taking
$A=0$ in Theorem 7.2 of \cite{bgm} shows that the Riemannian manifold 
\[(X_a,g):=((a,\infty)\times M, dt^2+\cosh(t)^2 h_M)\]
is a spin hyperbolic manifold of dimension $n+1$ for every
$a\in\rz\cup\{-\infty\}$.  
Setting $x:=e^{-t}$ near $t=\infty$, the metric $g$ becomes 
\[g=x^{-2}\left(dx^2+\frac{(1+x^2)^2}{4}h_M\right).\]
Theorem \ref{cor2} thus
shows that the Dirac operator on $X_a$ (or on any spin Riemannian
manifold containing $X_a$ as an open set) does not have $L^2$ eigenspinors
of real eigenvalue. This result was previously known when $M$
is compact. Interesting non-compact cases are obtained when $M$ is complete,
or when $M$ is compact with conical singularities with small angles 
\cite{chou}.

\subsection*{Rotationally symmetric Riemannian manifolds}
It is proved in \cite{anghel} that on $\rz^n$ with a metric which written in 
polar coordinates has the form $ds^2=dr^2 + \psi(r)^2 d\theta ^2$,
there are no $L^2$ harmonic spinors. This metric is complete. 
By the change of variables
\[x(r):=\int_r^\infty \frac{ds}{\psi(s)},\]
Anghel's metric becomes a particular case of \eqref{met} with $\psi(r(x))$ 
in the r\^ole of $f(x)$ from Theorem \ref{th1}, provided that
$\int_1^\infty dr/\psi(r)<\infty$. 
The absence of harmonic spinors is guaranteed in this case
by \cite{lott} if the resulting
conformal factor is Lipschitz. By Theorem \ref{th1} we know, even without 
the Lipschitz hypothesis, not only that there cannot exist $L^2$ harmonic 
spinors, but also that there are no $L^2$ eigenspinors at all.

\subsection*{The $L^2$-index of the Dirac operator}
Let $D^+$ denote the chiral component of $D$, viewed as an unbounded operator
in $L^2$, acting on compactly supported smooth spinors on a spin
Riemannian manifold as in Theorem \ref{th1}. Denote by $\oD^+$ its closure.
The $L^2$-index is defined as
\[\index(\oD^+):=\dim\ker (\oD^+)-\dim\ker (\oD^+)^*\]
where $(\oD^+)^*$ is the adjoint of $\oD^+$ (the definition makes
sense whenever 
both kernels are finite-dimensional, even when $\oD^+$ is not Fredholm). Here
$\ker(\oD^+)^*$ is precisely the distributional null-space of $D^-$ inside
$L^2$, while $\ker(\oD^+)$ is a subspace of the distributional null-space of 
$D^+$ inside $L^2$. Both these spaces vanish by Theorem \ref{th1}, so
in particular it follows that $\index(\oD^+)=0$.

\bibliographystyle{amsplain}

\end{document}